%% file: CatTopology.tex
\documentclass{llncs}
\usepackage[german,english]{babel}
\usepackage{times,theorem}
\usepackage{latexsym,color,graphics}
\usepackage{diagrams}
\diagramstyle[tight,centredisplay%
,dpi=600
,UglyObsolete
,heads=LaTeX%
]

\usepackage{epsfig}
\usepackage{amssymb}
\usepackage{amsmath}
\usepackage{amsfonts}
\usepackage{xspace}
\usepackage{graphicx}
\input{macros-general}

\title{Category Theory: Symmetry Group of  Comma-propagation Transformations}
\author{Zoran Majki\'c}
\authorrunning{Zoran Majki\'c}
\institute{ISRST, Tallahassee, FL, USA\\
\email{majk.1234@yahoo.com}\\ http://zoranmajkic.webs.com/}

\newtheorem{theo}{Theorem}
\newtheorem{propo}{Proposition}
\newtheorem{coro}{Corollary}

\begin{document}

\maketitle

\begin{abstract}
In general, all constructions of algebraic topology are functorial; the notions of category, functor and natural transformation originated here. The arrow categories are more simple forms of the \emph{comma} categories and were introduced by Lawvere  in the context of the interdefinability of the universal concepts of category theory. The basic idea is the elevation of arrows of one category $\textbf{C}$ to objects in another.
Given a category (as a "geometric object") $\textbf{C}$ we can consider its properties (the universal categorial commutative diagrams) preserved under  actions  of a comma-propagation operation $\{\}$ in the infinite hierarchy of its arrow-categories (n-dimensional levels, such that  for any $n\geq 1$,  $\textbf{C}_{n+1} = {\textbf{C}_n}$, with $\textbf{C}_1 =\textbf{C}$)  and on the functors (and their natural transformations) between such n-dimensional levels, which is a phenomena of a general categorial symmetry under a categorial-symmetry group $CS(\mathbb{Z})$ of all comma-propagation transformations.
\end{abstract}

\section{Introduction to Categorial Topology of n-dimensional Levels}
Most related philosophical questions deal with specific symmetries, objectivity,
interpreting limits on physical theories, classification, and laws of nature \cite{Shaw18,Eins36,Wign60,Koll10}. The recent history of the philosophy of mathematics is largely focused on grasping and defining the nature and essence of mathematics and its objects \cite{Korz94,Mein904,Bren874}.
Attempts to do this include explicating versions of: mathematics is just logic, mathematics is just structure, mathematics is a meaningless game, mathematics is a creation of the mind, mathematics is a useful fiction, etc.
Questions about how we can know numbers, shapes, and rules, and how it is even possible to prove necessary arithmetical, algebraic, and geometric relations; about
whether mathematics is an internal cognitive structure (that
evolved), a learned one and/or one existing in an independent realm; and about whether we can build a machine or a computer simulation that emulates our knowledge
structures and conscious experience, are quite old \cite{Kahn11,Jack83,Marc18,Zalt88}.

In Categorial Topology, given a category (as a "geometric object") \cite{Cart01,GrDi60,Vacu21,DeRi15,Majk24G}  we can consider its properties preserved under continuous action (a "deformation") of a comma-propagation operation.
However, the Metacategory space, valid for all  categories,  can not be defined by using well-know Grothendeick's approach with discrete ringed spaces as demonstrated in \cite{Majk24G}. However we can consider any category $\textbf{C}$ as an abstract geometric object,
 that is, a discrete space where the points  are the objects of this category and arrows between objects as the paths. Based on this approach, we  can define the Cat-vector space $V$ valid for all categories with noncommutative (and partial) addition operation for the vectors, and their inner product. For the categories where we define the norm ("length") of the vectors in $V$ we can define also the outer (wedge) product of the vectors in $V$ and we show that such Cat-algebra satisfies two fundamental properties of the Clifford geometric algebra \cite{Majk24G}.

So, in this paper we extend the work about Category symmetries provided in the  book \cite{Majk23s} by providing the category-symmetry group $CS(\mathbb{Z})$ and monoid $CS(\mathbb{N})$, and Peano-like axiom system for n-dimensional levels.
 We will  consider in this paper only the global categorial symmetries valid \emph{for all} categories.


 In general, all constructions of algebraic topology are functorial; the notions of category, functor and natural transformation originated here. The arrow categories are more simple forms of the \emph{comma} categories and were introduced by Lawvere \cite{Lawve63} in the context of the interdefinability of the universal concepts of category theory. The basic idea is the elevation of arrows of one category $\textbf{C}$ to objects in another.

It is well known that, given a base category $\textbf{C}$, we can represent its morphisms as objects by using a derived \emph{arrow} category $\textbf{C}\downarrow \textbf{C}$ (a special case of the comma category), such that for a given arrow $f:A \rightarrow B$ between two objects $A,B \in Ob_C$, we have the object denoted by $\langle A,B,f\rangle  \in Ob_{C\downarrow C}$ in this arrow category.
 What is important for functorial semantics of derivation of the objects in a given category $\textbf{C}$ from its arrows, is that it \emph{needs} to use necessarily the arrow categories $\textbf{C}\downarrow\textbf{C}$ where the arrows $g:a\rightarrow b$ of $\textbf{C}$ are encapsulated as the objects in $\textbf{C}\downarrow\textbf{C}$ denoted by $J(g)$, so the functorial representation has to be given by the functors $F:\textbf{C}\downarrow\textbf{C}\rightarrow \textbf{C}$.

 With this, we introduce the infinite hierarchy of the arrow categories, for a given base category  $\textbf{C}$, as follows.
The notion of hierarchy of common arrow categories for a given base category $\textbf{C}$ is inductively defined for $n\geq 1$ by
      \begin{equation}\label{eq:n-levels0}
      (\textbf{C}\downarrow\textbf{C})^1 =_{def} \textbf{C}\downarrow\textbf{C},~~~~~~  (\textbf{C}\downarrow\textbf{C})^{n+1} =_{def} (\textbf{C}\downarrow\textbf{C})^n \downarrow (\textbf{C}\downarrow\textbf{C})^n
      \end{equation}
and, derived from it, the \emph{n-dimensional levels} based \emph{only} on the base category $\textbf{C}$, instead on the its first arrow category $(\textbf{C}\downarrow\textbf{C})$ used in (\ref{eq:n-levels0}) for $n\geq 1$, denoted by
      $\textbf{C}_{n+1} =_{def} (\textbf{C}\downarrow\textbf{C})^{n}$
That is, the n-dimensional levels are defined inductively for $n\geq 1$ by
\begin{equation}\label{eq:n-levels0s2}
\textbf{C}_1 =_{def}  \textbf{C},  ~~~~~~~~\textbf{C}_{n+1}=_{def} \textbf{C}_{n}\downarrow\textbf{C}_{n}
\end{equation}
that is, inductively,  (n+1)-dimensional level is the arrow category obtained from the n-dimensional level, where the arrows of the n-dimensional level are composed by $2^{n-1}$ arrows of the base category $\textbf{C}$ (considered as the first 1-dimensional level as well which represents the leafs of the \emph{binary tree} of this n-dimensional syntax structure).\\
\textbf{Remark}: Note that such an infinitary hierarchy of n-dimensional levels there exist for any category. So it holds also for the minimal (non empty) discrete category $\textbf{1}$ composed by one unique object and unique identity arrow of this object. This infinitary n-dimensional levels of the base category $\textbf{1}$ there exist just on the compositionality of the identity arrow with itself, and show to us how so atomic (minimal) category generates the \emph{infinity} of the arrow categories derived from $\textbf{1}$. The existence of the identity arrows for each object in \emph{any} other finite or infinite category guarantees the existence of the infinitary n-dimensional hierarchy for each category.
\\$\square$\\
So, this phenomena is a \emph{foundational} and characteristic (distinguishable) property of each category and hence \emph{of the Category Theory}, similar to the foundational role of sets in the traditional set-based Mathematics or to the primitive built-in binary identity predicate in the First-order Logic. In what follows it will be used for the new multi-dimensional categorial definition of the set of natural numbers $\mathbb{N}$, different from the standard von Neumann set-based flat bidimensional definition\footnote{In defining the natural numbers, von Neumann begins by examining the most fundamental set, the empty set denoted by $\{\}$ and create the following inductive patern for the set-based definition of natural numbers:
\begin{itemize}
  \item number zero, $0 = \{\}$ has zero lements;
  \item number one, $1 = \{\{\}\}$ has one element(the empty set);
  \item number two, $2 = \{\{\},\{\{\}\}\}$ has two elements (the empty set and the set containing the empty set);
\end{itemize}
and this process would continue infinitely until all the natural numbers have been defined. Note that in this set-based definition of natural numbers is introduced the 'member of' relation $\in$, for example, $0 \in 1$, $0\in 2$ and $1\in 2$, etc...
}.

Thus, we can define inductively the set of natural number in a categorial way based on the base category  $\textbf{C} = \textbf{1}$ (and hence in a non set-based way, as in von Neumann inductive definition) in the following way: the number zero is defined by empty category $\textbf{0}$, while other are defined inductively from the base category $\textbf{1}$  (with single object denoted by $\bullet$ and only identity arrow $id$ of it), that is,

$0 = \textbf{0}$

$1 = \textbf{1}$

$2 = \textbf{1}\downarrow\textbf{1}$

$3 = (\textbf{1}\downarrow\textbf{1})\downarrow(\textbf{1}\downarrow\textbf{1})$

$4 = ((\textbf{1}\downarrow\textbf{1})\downarrow(\textbf{1}\downarrow\textbf{1}))\downarrow ((\textbf{1}\downarrow\textbf{1})\downarrow(\textbf{1}\downarrow\textbf{1}))$

etc...\\
That this inductive definition is not set-based but n-dimensional geometric definition is easy to verify by the consideration of what geometrically represent the arrows in the arrow categories which define  the natural numbers. To render more easy this presentation, in the place of single identity arrow $id:\bullet\rightarrow\bullet$ the labeled arrows like, $f_i,g_i,h_i,k_i,...$. Heaving this in mind, that each labeled arrow is equal to the identity arrow, we obtain the following multi-dimensional representation of natural numbers:

1. The number zero equal to empty category $\textbf{0}$ geometrically is zero-dimensional (without any dimension);

2. The number 1, equal to category $\textbf{1}$, where each (that is, unique) arrow is \emph{one-dimensional} (geometrically a linear line segment)
$$\bullet \rTo^{id} \bullet$$

3. The number 2, equal to arrow category $\textbf{1}\downarrow\textbf{1}$, where each (again unique) arrow $(h_1;k_1):J(f_1)\rightarrow J(g_1)$ where $J(f_1) = J(g_1)$ are equal to $J(id)$ (the unique object  in $\textbf{1}\downarrow\textbf{1}$), which in base category $\textbf{1}$ must satisfy the equation $k_1\circ f_1 = g_1\circ h_1$, that is  commutative diagram geometrically represented in a \emph{bidimensional} surface
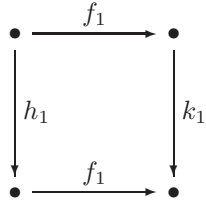
\begin{diagram}
 \bullet         &\rTo^{f_1}          & \bullet    \\
            \dTo_{h_1}&  &\dTo_{k_1}     \\
   \bullet        &    \rTo^{f_1}             &   \bullet
\end{diagram}

3. The number 3, equal to arrow category $(\textbf{1}\downarrow\textbf{1})\downarrow(\textbf{1}\downarrow\textbf{1})$, where each (again unique) arrow $((l_1;l_2);(l_3;l_4)):J(h_1;k_1)\rightarrow J(h_2;k_2)$ where $J(h_1;k_1) = J(h_2;k_2)$ are equal to $J(id;id)$ (the unique object  in $(\textbf{1}\downarrow\textbf{1})\downarrow(\textbf{1}\downarrow\textbf{1})$), which in base category $\textbf{1}$ must satisfy the six equations (of each external surface of the cube in next diagram, that is  commutative diagrams geometrically represented in the following  \emph{3-dimensional} cube
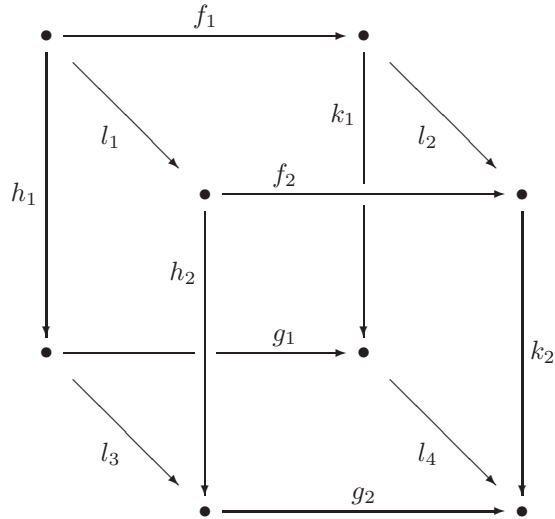
\begin{diagram}
\bullet & & \rTo^{f_1} & & \bullet & & \\
& \rdTo_{l_1} & & & \vLine^{k_1} & \rdTo_{l_2} & \\
\dTo^{h_1} & & \bullet & \rTo^{f_2} & \HonV & & \bullet \\
& & \dTo^{h_2} & & \dTo & & \\
\bullet & \hLine & \VonH & \rTo^{g_1} & \bullet & & \dTo_{k_2} \\
& \rdTo_{l_3} & & & & \rdTo_{l_4} & \\
& & \bullet & & \rTo^{g_2} & & \bullet \\
\end{diagram}
representing six \emph{equations}, the first two represented by the objects $J(h_1;k_1)$ and $J(h_2;k_2)$ and the next four derived from the arrow $((l_1;l_2);(l_3;l_4))$ from the first into second object:

$g_1\circ h_1 = k_1\circ f_1$

$g_2\circ h_2 = k_2\circ f_2$

$f_2\circ l_1 = l_2\circ f_1$

$h_2\circ l_1 = l_3\circ h_1$

$k_2\circ l_2 = l_4\circ k_1$

$g_2\circ l_3 = l_4\circ g_1$\\
etc...

Note that this definition is multi-dimensional, and not set-based, because the commutative diagrams are not the sets but the complex composition of arrows. More over, differently from the von Neumann set-based definition, we do not introduce the 'member of' relation, and in fact $0 \in 1$, and $1\in 2$ does not hold for this categorial multidimensional commutative-diagram's definition of natural numbers.

The commutative diagrams above distinguish clearly the Category Theory from the set-based theory but also from the pure geometric theory: it is true that geometrically these commutative diagrams have \emph{n-dimensional forms}, but they are also commutative and hence represents the \emph{algebraic equations}. I remarked previously that this pattern for definition of natural numbers can be defined for \emph{each} finite or infinite category $\textbf{C}$, so that in that case we have not only the commutative diagrams in such categories composed by only identity arrows (the trivial case), but also by the non-identity arrows! That is, the algebraic equations of the commutative diagrams in the pattern above are not trivial only but can be significative for a model theory represented by such a category  $\textbf{C}$ different from the minimal non-empty category $\textbf{1}$.

From my point of view, this pattern above demonstrates clearly how the Category Theory  as a particular synthesis of geometric and algebraic equational theory, is distinguished (its own) mathematical theory different both from the pure geometry and from the set-based mathematics.

We have provided not only a non set-theoretical definition of natural numbers, but have shown how the commutative diagrams in any base category $\textbf{C}$ are the sources for a generation of the multi-dimensional cubes, and this basic categorial 'spatial' (geometric) topology inside this base category is classified by its n-dimensional arrow categories: in each n-dimensional arrow category derived from a base category $\textbf{C}$, we have only the arrows that represent the n-dimensional commutative cubes of the base category $\textbf{C}$.  So, topologically a n-dimensional level represents a subcategory of $\textbf{C}$, composed by the n-dimensional commutative cubes.

This categorial partition of a base category is categorial analogous to the partition of a given set $S$ by its powerset $\P(S)$ whose elements are all subsets of $S$.

This presentation was a simple example of how we are using n-dimensional levels (of the arrow categories) in order to obtain mathematical frame for definition of the categorial symmetry. In what follows we will use the following standard comma projection functors for all n-dimensional levels of a given base category $\textbf{C}$, and $n\geq 1$:
\begin{itemize}
  \item The first comma projection, $F_{st} = (F_{st}^0,F_{st}^1):\textbf{C}_{n+1}\rightarrow \textbf{C}_{n}$,
      such that:\\
       for each object $\langle a,b,f\rangle$ in $\textbf{C}_{n+1}$,  $F_{st}(\langle a,b,f\rangle) = F_{st}^0(\langle a,b,f\rangle) = a$;\\
       for each arrow $(h_1,h_2)$ in $\textbf{C}_{n+1}$,  $F_{st}(h_1,h_2) = F_{st}^1(h_1,h_2) = h_1$.
  \item The second  projection, $S_{nd} = (S_{nd}^0,S_{nd}^1):\textbf{C}_{n+1}\rightarrow \textbf{C}_{n}$,
      such that:\\
       for each object $\langle a,b,f\rangle$ in $\textbf{C}_{n+1}$,  $S_{nd}(\langle a,b,f\rangle) = S_{nd}^0(\langle a,b,f\rangle) = b$;\\
       for each arrow $(h_1,h_2)$ in $\textbf{C}_{n+1}$,  $S_{nd}(h_1,h_2) = S_{nd}^1(h_1,h_2) = h_2$.
   \item The natural transformation  $\psi:F_{st} \rTo^\centerdot S_{nd}$.\\
              Note that a natural transformation $\psi$ associate to every object $X$ an arrow $\psi_X:F_{st}^0(X) \rightarrow S_{nd}^0(X)$, and this mapping can be represented by the function denoted by
       \begin{equation} \label{projectionTrnasf}
       J^{-1}:Ob_{\textbf{C}_{n+1}} \rightarrow Mor_{\textbf{C}_{n}}
       \end{equation}
For example, we have that $J^{-1}(\langle a,b,f\rangle) = \psi_{\langle a,b,f\rangle} =f$.
 \item We introduce the \emph{encapsulation} operator operator $J:Mor_{\textbf{C}_{n}}\rightarrow Ob_{\textbf{C}_{n+1}}$, as operator inverse to the operator $J^{-1}$, such that for each arrow $f:a\rightarrow b$ in $\textbf{C}_{n}$, we obtain       $J(f) = \langle a,b,f\rangle$.
\end{itemize}
The main concept in the n-dimensional levels theory of arrow categories is the fundamental coma-propagation process, based on the coherent transportation of the categories, functors and natural transformations for a given n-dimensional into next (n+1)-dimensional level of this stratified theory.

\section{Introduction to Comma Lifting and Comma-propagation in the n-dimensional  Levels Hierarchy \label{sec:nDimLev}}

We will use the incrementing symbol $~\widehat{}~~$ to denote such "coma lifting" of the categories, functors between them and natural transformations between such functors, as follows:
\begin{enumerate}
  \item For a category $\textbf{D}$, its comma lifting is $\widehat{\textbf{D}} =_{def} \textbf{D}\downarrow \textbf{D}$;
  \item For a covariant functor $F:\textbf{B}\rightarrow \textbf{D}$, its comma lifting functor is $\widehat{F}:\widehat{\textbf{B}}\rightarrow \widehat{\textbf{D}}$;
  \item For a natural transformation $\tau:F\rTo^\centerdot G$, its coma lifting natural transformation is $\widehat{\tau}:\widehat{F}\rTo^\centerdot \widehat{G}$.
\end{enumerate}
Note that the comma lifting represents the transfer of the semantics from a given theory into its metatheory.
This fact is easy to understand if we consider that each arrow $f:a\rightarrow b$ in the 'theory' category  $\textbf{B}$, equivalent to the commutative diagram $f\circ id_a = id_b\circ f$ where $id_a$ and $id_b$ are the identity arrows of the domain and codomain objects of this arrow $f$,   corresponds exactly to the arrow $(f;f):J(id_a) \rightarrow J(id_b)$ in the comma lifted metatheory category $\widehat{\textbf{B}} =(\textbf{B}\downarrow \textbf{B})$, and by using this correspondence, for each commutative diagram in the theory category $\textbf{B}$ (expressing some property of this theory)  there exists exactly the same commutative diagram in the metatheory (comma lifted) category  $\widehat{\textbf{B}} =(\textbf{B}\downarrow \textbf{B})$, by substituting each arrow $f$ in this diagram in theory category $\textbf{B}$, by the arrow $(f;f)$ in the metatheory category $\widehat{\textbf{B}} =(\textbf{B}\downarrow \textbf{B})$. Consequently,  by this process of comma lifting we support the propagation of the properties of a given theory to all higher metateories of this base theory.
These facts explain why the n-levels category theory has a general importance in the mathematics.


The following figure show the comma lifting of the functor $F:\textbf{B}\rightarrow\textbf{D}$ (between two commutative diagrams in the center of this figure) into the functor $\widehat{F}:\widehat{\textbf{B}}\rightarrow\widehat{\textbf{D}}$ (from left hand side arrow into right hand side arrow):
\begin{diagram}
J(f)&a         &\rTo^{f}          & b& F(a)         &\rTo^{F(f)}          & F(b)    & &&  \widehat{F}^0(J(f))= J(F^1(f))\\
\dTo^{(h;k)}&\Leftrightarrow \dTo_{h}       & &               \dTo^{k}&F~ \mapsto    \dTo_{F(h)}       & &               \dTo^{F(k)} &\Leftrightarrow &&  \dTo^{\widehat{F}^1(h;k)}_{=(F^1(h);F^1(k))} \\
J(g)&c  &  \rTo^{g} &   d& F(c)  &  \rTo^{F(g)}&   F(d) & && \widehat{F}^0(J(g)) = J(F^1(g)) \\
&& \textbf{B} && \rTo^F&  \textbf{D}  && && \\
\widehat{\textbf{B}} =(\textbf{B}\downarrow \textbf{B})&&  && \rTo^{\widehat{F}}&    && && \widehat{\textbf{D}}=(\textbf{D}\downarrow \textbf{D})
\end{diagram}
The comma lifting $\widehat{\tau}$ of the natural transformation $\tau:F\rTo^\centerdot G$ can be represented by the following figure  (in next we will use the representation of natural transformations as \emph{functions} from the object of the domain category  into the arrows of the codomain category of functors, that is, instead of its components $\tau_a$ for an object $a$, by $\tau(a)$):
\begin{diagram}
J(f)&&&\widehat{F}^0(J(f))= J(F^1(f))        &&&\rTo^{\widehat{\tau}_{J(f)}= \widehat{\tau}(J(f))}          &&& \widehat{G}^0(J(f))= J(G^1(f))\\
\dTo^{(h;k)}&&& \dTo^{\widehat{F}^1(h;k)}_{=(F^1(h);F^1(k))}       & &&&&&               \dTo^{\widehat{G}^1(h;k)}_{=(G^1(h);G^1(k))}\\
J(g)&&&\widehat{F}^0(J(g))= J(F^1(g))  &&&  \rTo^{\widehat{\tau}_{J(g)}= \widehat{\tau}(J(g))} &&&   \widehat{G}^0(J(g))= J(G^1(g))\\
in ~\widehat{\textbf{B}} =(\textbf{B}\downarrow \textbf{B}) &&& &&& in ~ \widehat{\textbf{D}} =(\textbf{D}\downarrow \textbf{D}) &&&
\end{diagram}
where the components of this comma lifted natural transformation are $\widehat{\tau}_{J(f)}= \widehat{\tau}(J(f)) =_{def} (\tau F_{st}^0; \tau S_{nd}^0)(J(f))$, that is,  in functional representation of natural transformations, we have that the comma lifted transformation is the following couple of functions:

$\widehat{\tau} =_{def} (\tau F_{st}^0; \tau S_{nd}^0):Ob_{\textbf{B}\downarrow\textbf{B}} \rightarrow Mor_{\textbf{D}\downarrow\textbf{D}}$.
\begin{example} \label{exam:nat-transf-lift}
Note that the commutative diagram in the figure above is just a commutative cube in the category $\textbf{D}$ composed by the following six equations:

1. Functorial consequences of the same base'theory' equation $k\circ f = g\circ h$: the $F(k)\circ F(f) = F(g)\circ F(h)$  and $G(k)\circ G(f) = G(g)\circ G(h)$ on the two opposite sides of the cube;

2. Consequences of the two arrows components of the 'theory' natural transformation $\tau$ for the arrows $f:a\rightarrow b$ and $g:c\rightarrow d$ in $\textbf{B}$: the $G(f)\circ \tau(F_{st}^0J(f)) = G(f)\circ \tau(a)= \tau(b)\circ F(f) = \tau(S_{nd}^0J(f))\circ F(f)$  obtained from the arrow $f:a\rightarrow b$, $G(g)\circ \tau(F_{st}^0J(g)) = G(g)\circ \tau(c)= \tau(d)\circ F(g) = \tau(S_{nd}^0J(g))\circ F(g)$, on the two opposite sides of the cube;

3. Consequences of the two arrows components of the 'theory' natural transformation $\tau$ for the arrows $h:a\rightarrow c$ and $k:b\rightarrow d$ in $\textbf{B}$:: the $G(h)\circ \tau(F_{st}^0J(f)) = G(h)\circ \tau(a) = \tau(c)\circ F(h)= \tau(F_{st}^0J(g)) \circ F(h)$ and $G(k)\circ \tau(S_{nd}^0J(f)) = G(k)\circ \tau(b) = \tau(d)\circ F(k)= \tau(S_{nd}^0J(g)) \circ F(k)$,  on the two opposite sides of the cube.
\end{example}
\textbf{Remark}: Thus, all equations are derived by natural transformation and two functors, from the unique base 'theory' equation (commutative diagram) $k\circ f = g\circ h$ in $\textbf{B}$, so that the 'set of correlated equations' is this result of the action of the natural transformation $\tau$ and its two functors $F$ and $G$.

So, as shown in the last figure above, from the fact that the (unique) arrow $(h;k):J(f) \rightarrow J(g)$ in the comma lifted category (metatheory)  $\widehat{\textbf{B}}  = (\textbf{B}\downarrow\textbf{B})$ is just an equation (commutative diagram in 'theory' $\textbf{B}$, the commutative diagram obtained from this arrow by action of the comma lifted natural transformation $\widehat{\tau}$  in effect represents these six correlated equations (commutative diagrams in denotational semantics category $\textbf{D}$) in a compact representation by unique commutative diagram in the comma lifted $\widehat{\textbf{D}}  = (\textbf{D}\downarrow\textbf{D})$.
\\$\square$
\begin{propo}\label{propA.1} \cite{Majk23s}
The operator of "comma lifting" is defined for the categories by $\widehat{\textbf{B}}  = (\textbf{B}\downarrow\textbf{B})$. The following functors and natural transformation there exist as an comma lifting consequence in the higer n-dimensional levels:

 For each functor $F = (F^1,F^2):\textbf{B}\rightarrow\textbf{D}$, there exists the following comma lifted functor
      \begin{equation}\label{eq:commaFunctor}
      \widehat{F} =_{def} (JF^1\psi,(F^1F_{st}^1;F^1S_{nd}^1)):\widehat{\textbf{B}}\rightarrow\widehat{\textbf{D}}
      \end{equation}
      and for each natural transformation $\tau:F\rTo^\centerdot G$ there exists the following comma lifted natural transformation
      \begin{equation}\label{eq:commaNatTr}
      \widehat{\tau} =_{def} (\tau F_{st}^0;\tau S_{nd}^0):\widehat{F}\rTo^\centerdot \widehat{G}
      \end{equation}
      and the following functor denominated "nat-functor"
      \begin{equation}\label{eq:NatTrFunctor}
      N_\tau =_{def}(J\tau, \widehat{\psi}J\widehat{\tau}J):\textbf{B}\rightarrow(\textbf{D}\downarrow\textbf{D})
      \end{equation}
      where the natural transformation $\widehat{\psi}:\widehat{F}_{st}\rTo^\centerdot \widehat{S}_{nd}$ is the comma lifted \emph{fundamental} natural transformation $\psi$.
  \end{propo}
Note that we denominated the functor $N_\tau$ by "nat-functor" because it is a direct consequence of the natural transformation $\tau$. We have seen that the duality functor is a particular case of such general nat-functor, so let us consider the other examples of this nat-functor:
\begin{example} \label{exam:diag-functor}
Let us consider the following  case of the general nat-functor $N_\tau$ in (\ref{eq:NatTrFunctor}, when $\textbf{B} = \textbf{D}$, the nat-functor $N_{\tau_I}$  in (\ref{eq:NatTrFunctor}), obtained from the identity natural transformation $\tau_I:Id_\textbf{C}\rTo^\centerdot Id_\textbf{C}$  obtained from the identity functor $F = G = Id_\textbf{C}:\textbf{C}\rightarrow\textbf{C}$. \\It is easy to show that this nat-functor is equal to the well known \emph{diagonal} functor (we will only substitute a category $\textbf{B}$ by $\textbf{C}$, in order to obtain the same original syntax used for n-dimensional levels), so from (\ref{eq:NatTrFunctor}),
\begin{equation} \label{eq:diagFunc}
\blacktriangle = (J\tau_I, \widehat{\psi}J\widehat{\tau}_I J):\textbf{C}\rightarrow(\textbf{C}\downarrow\textbf{C})
\end{equation}
such that for any object $a \in \textbf{C}$, it holds that $\blacktriangle^0(a) = J\tau_I(a) = J(id_a) = \langle a,a, id_a\rangle$, and for each arrow $f:a\rightarrow b$, from the derivation in the prof above we obtain $\blacktriangle^1(f) = \widehat{\psi}J\widehat{\tau}_I J(f) = (F^1;G^1)(f) = (Id_\textbf{C}^1;Id_\textbf{C}^1)(f) = (f;f)$.
\end{example}
$\square$\\
 The next theorem \cite{Majk23s} demonstrates one of the fundamental structural properties preserved by the comma lifting:
 \begin{theo} \label{th:A1} \cite{Majk23s}
 Given a functor $G:\textbf{D}\rightarrow \textbf{C}$, for each universal arrow $(d,g)$, with $d\in Ob_{\textbf{D}}$ and arrow $g:c\rightarrow G(d)$ in $\textbf{C}$,  from $c$ to functor $G$, presented by the figure\footnote{Where for each pair $(d'_i,f_i)$ there exists \emph{the unique} arrow $\underline{f}_i:d\rightarrow d'_i$ such that $f_i = G(\underline{f}_i)\circ g$.},
  \begin{equation} \label{fig:universal-arrow2}
\begin{diagram}
 c         &\rTo^{g}          & G(d)    && d\\
           & \rdTo_{f_i}&     \dTo_{G(\underline{f}_i)}&&\dDashto_{\underline{f}_i}     \\
           &                 &   G(d'_i)   &&d'_i
\end{diagram}
\end{equation}
there is the universal arrow $(J(id_d),(g,g))$ from $J(id_d)$ to comma lifted functor $\widehat{G}:(\textbf{D}\downarrow\textbf{D})\rightarrow (\textbf{C}\downarrow\textbf{C})$.
If $F:\textbf{C}\rightarrow \textbf{D}$ is left adjoint functor to $G$, with adjunction $(F,G,\varepsilon,\eta)$ with natural transformations, unit $\eta:Id_\textbf{C}\rTo^\centerdot GF$ and  counit $\varepsilon:FG\rTo^\centerdot Id_\textbf{D}$ then there exists the comma-lifted adjunction $(\widehat{F},\widehat{G},\widehat{\varepsilon},\widehat{\eta})$.
 \end{theo}
%
 %
 \textbf{Comma-propagation:}  The phenomena of the comma-propagation in the n-dimensional hierarchy, for the functors and their derived concepts as natural transformations and adjunctions, is based on the effects of the harmonics in the given spectra  (like the spectra of the sounds, waves, etc..). The existence of a functor between two basic categories (1-dimensional level) causes the existence of comma lifted functors between the arrow categories of higher n-dimensional levels.

 Comma-propagation is a \emph{general transformation} applied to n-dimensional levels, their functors and natural transformations. In what follows, the operation of the comma-propagation will be denoted by $\{\_\}$, and hence applied to a category $\textbf{C}$, functor $F$ and natural transformation $\eta$ will be denoted by $\{\textbf{C}\}$, $\{F\}$ and $\{\eta\}$ relatively. We will introduce this operation for the n-dimensional levels by the following definition:
\begin{definition} \label{def:comProp} \textsc{Comma-propagation transformation of n-dimensional levels}:\\
For any given n-dimensional level $\textbf{C}_n$, we define its comma-propagation by $\{\textbf{C}_n\} \triangleq \textbf{C}_{n+1}$, such that
\begin{enumerate}
  \item For each object $c$ in $\textbf{C}_n$, the comma-propagated object in $\{\textbf{C}_n\}$ is defined by
      \begin{equation}\label{eq:CP-obj}
      \{c\} \triangleq \blacktriangle^0(c) = J(id_c)
      \end{equation}
  \item For each arrow $f$ in $\textbf{C}_n$, the comma-propagated arrow in $\{\textbf{C}_n\}$ is defined by
      \begin{equation}\label{eq:CP-arr0w}
      \{f\} \triangleq \blacktriangle^1(f) = (f;f)
      \end{equation}
\end{enumerate}
\end{definition}
Let us show that each (n+1)-dimensional level $\textbf{C}_{n+1}= \textbf{C}_n\downarrow\textbf{C}_n =(\textbf{C}\downarrow\textbf{C})^{n}$, for $n\geq 1$, with $\textbf{C}_{1} = \textbf{C}$, of a given base category $\textbf{C}$, can be equivalently represented by the category of functors $\textbf{C}_n^{\textbf{J}}$ where $\textbf{C}_n$  is n-dimensional level and the small index category $\textbf{J}$ is equal to the preorder category $\textbf{2}$ with two objects $a_1= 0$ and $a_2=1$ and unique non-identity arrow $l_{12}:a_1\rightarrow a_2$ (representing partial order $0\leq 1$).

 In fact, for an arrow (natural transformation) $\eta:F\rTo^\centerdot F'$ in $\textbf{C}_n^{\textbf{2}}$ between two objects $F$ and $F'$ which are functors $F,F':\textbf{2}\rightarrow \textbf{C}_n$, we have the corresponding arrow $(\eta(a_1);\eta(a_2)):J(F(l_{12}))\rightarrow J(F'(l_{12}))$ in  $\textbf{C}_{n+1}=(\textbf{C}\downarrow\textbf{C})^{n}$, corresponding to the following commutative diagram in  $\textbf{C}_{n}$
 \begin{diagram}
  F(a_1)       &&\rTo^{\eta(a_1)}          && F'(a_1)\\
 \dTo_{F(l_{12})}     & &&&     \dTo_{F'(l_{12})}     \\
 F(a_2)          && \rTo^{\eta(a_2)}        &&   F'(a_2)
\end{diagram}
 Viceversa, for each arrow $(h;k):J(f)\rightarrow J(g)$ in $\textbf{C}_{n+1}$, we have the arrow (natural transformation)
 $\eta:F\rTo^\centerdot F'$ in $\textbf{C}_n^{\textbf{2}}$ with components $\eta(a_1) = h$ and $\eta(a_2) = k$, between functor $F$ defined by $F(l_{12}) \triangleq f$ and functor $F'$ defined by $F'(l_{12}) \triangleq g$. Thus, $\textbf{C}_{n+1}$ is equivalent to $\textbf{C}_n^{\textbf{2}}$.
 It is well know that small index categories $\textbf{J}$ with a category of functors $\textbf{C}^{\textbf{J}}$ are used to define the \emph{limits} in the base category $\textbf{C}$. With comma-propagation we can see how the limits in the base category $\textbf{C}_1 =\textbf{C}$ are inductively propagated in all higher n-dimensional levels $\textbf{C}_n$.

 In fact, we can generalize the specific (for arrow categories) arrow-diagonal functor $\blacktriangle:\textbf{D}\rightarrow \textbf{D}\downarrow \textbf{D}$ and the standard diagonal functor $\vartriangle:\textbf{D}\rightarrow \textbf{D}\times \textbf{D}$,
 into a more \emph{general} diagonal functor (ascending case) $\vartriangle:\textbf{D}\rightarrow \textbf{C}^{\textbf{J}}$, where $\textbf{C} = \textbf{D}^m \triangleq~\overbrace{\textbf{D}\times...\times \textbf{D}}^m$ is n-ary  product of categories for finite $m\geq1$, as follows:
 \begin{definition}\label{def:genDiag}
 We define the general diagonal functor $\vartriangle = (\vartriangle^0,\vartriangle^1):\textbf{D}\rightarrow \textbf{C}^{\textbf{J}}$,
 where for a finite $m\geq 1$, $\textbf{C} = \textbf{D}^m \triangleq~\overbrace{\textbf{D}\times...\times \textbf{D}}^m$ with $\textbf{D}^1 = \textbf{D}$,
 such that
 \begin{enumerate}
   \item  For each object $c$ in $\textbf{D}$, we define the constant functor $\vartriangle^0(c):\textbf{J}\rightarrow \textbf{C}$, such that for all indexed objects $a_j$ and arrows $l_{jk}:a_j \rightarrow a_k$ in $\textbf{J}$,
       \begin{equation} \label{eq:genDiag0}
       \vartriangle^0(c)(a_j) = (\overbrace{c,...,c}^m)~~~~~ \emph{and} ~~~~~ \vartriangle^0(c)(l_{jk}) = id_{(c,...,c)}
       \end{equation}
   \item For each arrow $g:c \rightarrow c'$ in in $\textbf{D}$, $\vartriangle^1(g)$ is a constant natural transformations, such that for each indexed object $a_j$ in  $\textbf{J}$, the arrow component in $\textbf{C}$ of this natural transformation is
       \begin{equation} \label{eq:genDiag1}
       \vartriangle^1(g)(a_j) = (\overbrace{g,...,g}^m)
       \end{equation}
 \end{enumerate}
 \end{definition}
 The following figure with commutative diagram demonstrates the relationship between the arrow-diagonal functor $\blacktriangle:\textbf{D}\rightarrow (\textbf{D}\downarrow\textbf{D})$ and this more general diagonal functor $\vartriangle:\textbf{D}\rightarrow\textbf{C}^{\textbf{J}}$, for $m=1$, $\textbf{C} = \textbf{D}^1=\textbf{D}$,
 \begin{equation} \label{fig:arrow-diagonal}
 \begin{diagram}
c&&F =\vartriangle^0(c) &&&  F(a_j)=c       &&\rTo^{F(l_{jk})}          & F(a_k)=c&&  J(id_c)\\
\dTo^{g}&& \dTo^{\eta}_{=\vartriangle^1(g)}&&&\dTo^{\eta(a_j)}_{=g}     && &     \dTo^{\eta(a_k)}_{=g}   &\Leftrightarrow&  \dTo_{(g;g)}\\
c'&&G =\vartriangle^0(c')&&& G(a_j) =c'         && \rTo^{G(l_{jk})}        &   G(a_k)=c'&&J(id_{c'})\\
 \textbf{D} &\rTo^{\vartriangle} &\textbf{D}^{\textbf{J}} &&&&& in ~\textbf{D} &&& \\
  \textbf{D} && &&&\rTo^{\blacktriangle}&&  &&& (\textbf{D}\downarrow\textbf{D})
\end{diagram}
\end{equation}
Note that in the case when  the small index category $\textbf{J} = \textbf{2}$ is composed by only two objects $0$ and $1$ and non-identity arrow $e:0 \rightarrow 1$, in the commutative diagram above (for $j =1$ and $k=2$) we have that $a_j = a_1= 0$, $a_k =a_2 = 1$ and arrow $l_{jk}= l_{12} = e$, so we have the exact correspondence  between the objects $\vartriangle^0(c)$ and $J(id_c)$ and arrows $\vartriangle^1(g)$ and $(g;g)$ in the categories $\textbf{D}^{\textbf{2}}$ and $(\textbf{D}\downarrow\textbf{D})$, respectively.

The following figure demonstrates the relationship between the standard diagonal functor $\vartriangle:\textbf{D}\rightarrow \textbf{D}\times\textbf{D}$ and this more general diagonal functor $\vartriangle:\textbf{D}\rightarrow\textbf{C}^{\textbf{J}}$, in the case when $m=2$, $\textbf{C} = \textbf{D}^2=\textbf{D}\times\textbf{D}$, and $\textbf{J} = \textbf{1}$ is the category of only one object $a_1$ and its identity arrow $id_{a_1}$,
\begin{equation} \label{fig:DiagGen}
 \begin{diagram}
c&&F =\vartriangle^0(c):\textbf{1}\rightarrow \textbf{C} &&&&  F(a_1)=(c,c)~~  &&\rTo^{F(id_{a_1})=id_{(c,c)}}&& ~~~F(a_1)=(c,c)\\
\dTo^{g}&& \dTo^{\eta}_{=\vartriangle^1(g)}&&&&\dTo^{\eta(a_1)}_{=(g,g)}     &&&&      \dTo^{\eta(a_1)}_{=(g,g)}   \\
c'&&G =\vartriangle^0(c'):\textbf{1}\rightarrow \textbf{C}&&&& G(a_1) =(c',c')~~~~         && \rTo^{G(id_{a_1})= id_{(c',c')}}&&  ~~~~ G(a_1)=(c',c')\\
 \textbf{D} &\rTo^{\vartriangle} &(\textbf{D}\times\textbf{D})^{\textbf{1}} &&&&&& in ~\textbf{C} = \textbf{D}\times\textbf{D}&&
 \end{diagram}
\end{equation}
 With this generalization, for the  comma-propagation, we will introduce two new specific functors, called \emph{modulators},  $L$ and $K$ as follows:
 \begin{definition} \label{def.K-L-iso}
 Let $\textbf{C}^{\textbf{J}}$ be a category of functors from some small index category $\textbf{J}$ with indexed objects $a_j$, $1\leq j\leq N$, and indexed arrows $l_{j,m}:a_j\rightarrow a_m$, into the category $\textbf{C}$. Then for each n-dimensional level $\textbf{C}_n$, $n\geq 1$, with $\textbf{C}_{n+1}=\textbf{C}_n\downarrow\textbf{C}_n$, there the following two functors: $L= (L^0,L^1):(\textbf{C}_n^{\textbf{J}}\downarrow \textbf{C}_n^{\textbf{J}})\rightarrow (\textbf{C}_n\downarrow\textbf{C}_n)^{\textbf{J}}$, and  $K:(\textbf{C}_n\downarrow\textbf{C}_n)^{\textbf{J}}\rightarrow (\textbf{C}_n^{\textbf{J}}\downarrow \textbf{C}_n^{\textbf{J}})$, defined by:
 \begin{enumerate}
   \item The object component $L^0$ of the functor $L:(\textbf{C}_n^{\textbf{J}}\downarrow \textbf{C}_n^{\textbf{J}})\rightarrow (\textbf{C}_n\downarrow\textbf{C}_n)^{\textbf{J}}$ is defined by:\\For every object $J(h) \in Ob_{(\textbf{C}_n^{\textbf{J}}\downarrow \textbf{C}_n^{\textbf{J}})}$, where $h:G\rTo^\centerdot H$ is a natural transformation between functors $G,H:\textbf{J}\rightarrow \textbf{C}_n$, we obtain the functor $F = L^0(J(h)):\textbf{J}\rightarrow (\textbf{C}_n\downarrow\textbf{C}_n)$, such that for each indexed object $a_j$ in $\textbf{J}$
 \begin{equation}\label{eq:L-object}
 F(a_j) \triangleq~ J(h(a_j))
 \end{equation}
 where $h(a_i):G(a_i)\rightarrow H(a_i)$ is the i-th arrow component of the natural transformation $h$.
  For each arrow $(g_1;g_2):J(h)\rightarrow J(k)$ in $(\textbf{C}_n^{\textbf{J}}\downarrow \textbf{C}_n^{\textbf{J}})$ which represents the commutative diagram of natural transformation in $\textbf{C}_n^{\textbf{J}}$, composed by vertical composition $g_2\bullet h =  k\bullet g_1$, such that for each indexed object $a_j$ in ${\textbf{J}}$, we have the commutative diagram in $\textbf{C}_n$, equivalent to vertical arrow in $(\textbf{C}_n\downarrow\textbf{C}_n)$,
 \begin{equation}\label{eq:vert-diag}
 \begin{diagram}
  G(a_j)       &\rTo^{h(a_j)}          & H(a_j)&&&& F(a_j) = J(h(a_j))\\
 \dTo^{g_1(a_j)}     & &     \dTo_{g_2(a_j)}   &&\Leftrightarrow& & \dTo^{\alpha(a_j)=}_{(g_1(a_j);g_2(a_j))}\\
 G'(a_j)          & \rTo^{k(a_j)}        &   H'(a_j)&&&&F'(a_j) = J(k(a_j))
\end{diagram}
\end{equation}
So, the arrow component $L^1$  is defined by a natural transformation $L^1(g_1;g_2) = \alpha:F\rTo^\centerdot F'$ between the functors $F = L^0(J(h))$ and $F' = L^0(J(k))$, such that for each indexed object $a_j$ in ${\textbf{J}}$,
 \begin{equation}\label{eq:L-object2}
 \alpha(a_j) \triangleq~ (g_1(a_j);g_2(a_j))
 \end{equation}
 is the vertical arrow in the right hand side of figure (\ref{eq:vert-diag}).
  \item The object component $K^0$ of the functor $K:(\textbf{C}_n\downarrow\textbf{C}_n)^{\textbf{J}}\rightarrow (\textbf{C}_n^{\textbf{J}}\downarrow \textbf{C}_n^{\textbf{J}})$ is defined by:\\For every object (functor) $F$ in $(\textbf{C}_n\downarrow\textbf{C}_n)^{\textbf{J}}$ we define $K^0(F) = J(h)$ where $h$ is an arrow in $\textbf{C}_n^{\textbf{J}}$ and hence a natural transformation between some functors $G,H:{\textbf{J}}\rightarrow \textbf{C}_n$, that must satisfy for each $a_j$-indexed  arrow component
  \begin{equation}\label{eq:L-object3}
        h(a_j) \triangleq~ \psi(F(a_j)
  \end{equation}
 i.e., $J(h(a_j) = F(a_j)$. Let, analogously, $K^0(F') = J(k)$  for an arrow (natural transformation) $k$ in $\textbf{C}_n^{\textbf{J}}$. \\For each arrow  $\alpha:F\rTo^\centerdot F'$ in $(\textbf{C}_n\downarrow\textbf{C}_n)^{\textbf{J}}$, which is a natural transformation, the arrow $K^1(\alpha):K^0(F)\rightarrow K^0(F')$ is an arrow in the arrow category $(\textbf{C}_n\downarrow\textbf{C}_n)^{\textbf{J}}$ and hence represented by a pair of arrows $(g_1,g_2) = K^1(\alpha)$, that is $K^1(\alpha):J(h) \rightarrow J(k)$, so that this arrow represents a commutative diagram of vertical composition of natural transformations $g_2\bullet h =  k\bullet g_1$.
 So, natural transformation $K^1(\alpha)$ is defined in the way that for each j-th object $a_j$ in ${\textbf{J}}$,
  \begin{equation}\label{eq:L-object4}
       (K^1(\alpha))(a_j) \triangleq~ (g_1(a_j);g_2(a_j))
  \end{equation}
 \end{enumerate}
 \end{definition}
 It is easy to verify that $L$ and $K$ are well defined, that is, they map  identity arrows  $(g_1,g_2)$ (when $h  = k$ and $g_1=dom(h)$, $g_2 = cod(h)$) into identity arrows.   The functorial property for the composition of arrows are satisfied by the commutativity of the extension of commutative diagrams in (\ref{eq:vert-diag}) for each indexed object $a_j$ in ${\textbf{J}}$. That is, for a composition of two arrows  $l\circ g$  in $(\textbf{C}_n^{\textbf{J}}\downarrow \textbf{C}_n^{\textbf{J}})$, with $g = (g_1;g_2)$ and $l=(l_1;l_2)$, let us show that $L(l\circ g) = L(l)\bullet L(g)$ is just the vertical composition of the natural transformations $L(l)$ and $L(g)$, by using their components for each object $a_j$ in $\textbf{J}$:

 $L(l\circ g)(a_j) =L((l_1;l_2)\circ (g_1;g_2))(a_j)$

 $= L(l_1\bullet g_1;l_2\bullet g_2)(a_j)$   because $l_1,l_2,g_1$ and $g_2$ are nat. transformations

 $= ((l_1\bullet g_1)(a_j);(l_2\bullet g_2)(a_j))$  from definition of $L$

 $= (l_1(a_j)\circ g_1(a_j); l_2(a_j)\circ g_2(a_j))$  from definition of vertical composition of nat.transf.

 $= (l_1(a_j);l_2(a_j))\circ (g_1(a_j);g_2(a_j))$

 $= L(l_1;l_2)(a_j) \circ L(g_1;g_2)(a_j)$.\\
 It is shown (in \cite{Majk23s}) that $K$ is right adjoint of $L$, and that they are isomorphisms.
 \begin{propo} \label{prop.A8} \cite{Majk23s}
 For each n-dimensional level $\textbf{C}_n$, $n\geq 1$, there the following isomorphism of the categories $L= (L^0,L^1):(\textbf{C}_n^{\textbf{J}}\downarrow \textbf{C}_n^{\textbf{J}})\rightarrow (\textbf{C}_n\downarrow\textbf{C}_n)^{\textbf{J}}$, and its right adjoint functor $K=(K^0,K^1):(\textbf{C}_n\downarrow\textbf{C}_n)^{\textbf{J}}\rightarrow (\textbf{C}_n^{\textbf{J}}\downarrow \textbf{C}_n^{\textbf{J}})$, with
 \begin{equation}\label{eq:L-object4}
       L\circ K =~ I_{(\textbf{C}_n\downarrow\textbf{C}_n)^{\textbf{J}}}~~~~ \emph{and}~~~~ K\circ L = I_{(\textbf{C}_n^{\textbf{J}}\downarrow \textbf{C}_n^{\textbf{J}})}
  \end{equation}
  where $I_{(\textbf{C}_n\downarrow\textbf{C})_n^{\textbf{J}}}$ and $I_{(\textbf{C}_n^{\textbf{J}}\downarrow \textbf{C}_n^{\textbf{J}})}$ are two identity functors, relatively. This adjunction is denoted by $(L,K,id_{LK},id_{KL})$ where the identity natural transformation $id_{KL}:I_{(\textbf{C}_n^{\textbf{J}}\downarrow \textbf{C}_n^{\textbf{J}})}\rightarrow KL$ is the unit and $id_{LK}:LK\rightarrow I_{(\textbf{C}_n\downarrow\textbf{C}_n)^{\textbf{J}}}$  is the counit of this adjunction.

  In the case when ${\textbf{J}} = \textbf{2}$ is the category of only two indexed objects $a_1$ and $a_2$ with unique non-identity arrow $l_{12}:a_1 \rightarrow a_2$, we obtain that

  $L = K:\textbf{C}_{n+2}\rightarrow \textbf{C}_{n+2}$\\
  that is, $K$ and $L$ are equal (the objects $\textbf{C}_n^{\textbf{2}}$ and $\textbf{C}_{n+1} = (\textbf{C}_n\downarrow\textbf{C}_n)$ are equal up to the isomorphisms  in the \textbf{Cat} category).
 \end{propo}
 So, we are able to define comma-propagation of functors and natural transformations as follows:
 \begin{definition}\label{def:commaPropagationFunctor}
 For a given small index category $\textbf{J}$, we define the following four cases of the comma-propagation  of a given base covariant functor $F$  (considered as 1-level) which for each n-dimensional level,   will be denoted by $F_n$ comma-propagated functor with the n-dimensional levels $\textbf{C}_n$ and $\textbf{D}_n$, as follows  for $n\geq 2$:
 \begin{enumerate}
   \item Basic case: when the base functor is $F:\textbf{D}\rightarrow \textbf{C}$, with $F_1 \triangleq~ F$, and let $F_{n-1}:\textbf{D}_{n-1}\rightarrow \textbf{C}_{n-1}$ be (n-1)-dimensional comma-propagated functor.  Then n-dimensional  comma-propagatet functor is its comma lifted\\
       $F_n = \{F_{n-1}\} \triangleq~ \widehat{F}_{n-1}:\textbf{D}_n\rightarrow \textbf{C}_n$.
   \item Descending case: when the base functor is $F:\textbf{D}^{\textbf{J}}\rightarrow \textbf{C}$, with $F_1 \triangleq~ F$, and let $F_{n-1}:\textbf{D}_{n-1}^{\textbf{J}}\rightarrow \textbf{C}_{n-1}$ be (n-1)-dimensional comma-propagated functor.  Then, for $K:\textbf{D}_{n}^{\textbf{J}}\rightarrow (\textbf{D}_{n-1}^{\textbf{J}}\downarrow\textbf{D}_{n-1}^{\textbf{J}})$,  n-dimensional comma-propagated functor is \\
       $F_n = \{F_{n-1}\} \triangleq~ \widehat{F}_{n-1} \circ K:\textbf{D}_{n}^{\textbf{J}}\rightarrow \textbf{C}_n$.
   \item Ascending case: when the base functor is $F:\textbf{D}\rightarrow \textbf{C}^{\textbf{J}}$, with $F_1 \triangleq~ F$, and let $F_{n-1}:\textbf{D}_{n-1}\rightarrow \textbf{C}_{n-1}^{\textbf{J}}$ be (n-1)-dimensional comma-propagated functor.  Then, for $L:(\textbf{C}_{n-1}^{\textbf{J}}\downarrow\textbf{C}_{n-1}^{\textbf{J}}) \rightarrow \textbf{C}_{n}^{\textbf{J}}$, n-dimensional  comma-propagated functor is \\
       $F_n = \{F_{n-1}\} \triangleq~ L\circ\widehat{F}_{n-1}:\textbf{D}_{n}\rightarrow \textbf{C}_{n}^{\textbf{J}}$.
   \item Balanced case: when the base functor is $F:(\textbf{D}\downarrow\textbf{D})\rightarrow (\textbf{C}\downarrow\textbf{C})$, with $F_1 \triangleq~ F$, and let $F_{n-1}:\textbf{D}_{n-1}^{\textbf{J}}\rightarrow \textbf{C}_{n-1}^{\mathfrak{J}}$ be (n-1)-dimensional comma-propagated functor.  Then, for $K:\textbf{D}_{n}^{\textbf{J}}\rightarrow (\textbf{D}_{n-1}^{\textbf{J}}\downarrow\textbf{D}_{n-1}^{\textbf{J}})$ and $L:(\textbf{C}_{n-1}^{\mathfrak{J}}\downarrow\textbf{C}_{n-1}^{\mathfrak{J}}) \rightarrow \textbf{C}_{n}^{\mathfrak{J}}$, n-dimensional  comma-propagated functor is \\
       $F_n = \{F_{n-1}\} \triangleq~ L\circ\widehat{F}_{n-1}\circ K:\textbf{D}_{n}^{\textbf{J}}\rightarrow \textbf{C}_{n}^{\mathfrak{J}}$\\
       where  ${\mathfrak{J}}$ is another small index category.
 \end{enumerate}
 The functors $K$ and $L$ in Definition \ref{def.K-L-iso} are initial/final \emph{modulators} relatively: $K$ is used to pass from the category of functors into an arrow category, from the fact that the comma lifted functor works only for the hierarchy of arrow categories, while $L$ is used to pass from an arrow category (the result of comma lifting) into final category of functors.
  \end{definition}
Consequently, in the basic case, the comma-propagated functor $F_n$ is just the n-times iteration of the comma lifting of a functor defined by (\ref{eq:commaFunctor}), $\widehat{G} \triangleq~ (JG^1\psi,(G^1F_{st}^1;G^1S_{nd}^1))$, exactly as we obtained the comma-induced functor $F_n$. For the remaining three cases, modulated by the functors $K$ and $L$ in Definition \ref{def.K-L-iso}, the comma-propagation reduces to the comma-induction (defined in \cite{Majk23s})in the particular cases when the small index categories $\textbf{J}$ and $\mathfrak{J}$ are equal to the one-arrow index category $\textbf{2}$.

Let us show in an example how the comma-propagation of functors preserves their original properties in the higher n-dimensional levels:
\begin{example} \label{exam:DiagFunct}
Let us consider the general diagonal functor (which is an ascending case) $\vartriangle:\textbf{D}\rightarrow \textbf{C}^{\textbf{J}}$ in Definition \ref{def:genDiag}, when $\textbf{C}=\textbf{D}$ with $m=1$, that is, $\vartriangle:\textbf{C}\rightarrow \textbf{C}^{\textbf{J}}$, and suppose that for $n\geq 1$ for comma-propagated diagonal functor  $\vartriangle_n:\textbf{C}_n\rightarrow \textbf{C}_n^{\textbf{J}}$ it holds the diagonal property, with $\vartriangle_n =\vartriangle$. That is, for each object $c \in Ob_{\textbf{C}_n}$, the object (a functor) $\vartriangle^0_n(c):\textbf{J}\rightarrow \textbf{C}_n$ satisfies the diagonal property in point 1 of Definition \ref{def:genDiag}, such that for all indexed objects $a_i$ in $\textbf{J}$,

$\vartriangle^0_n(c)(a_j) = c$\\
and that for each arrow $g \in Mor_{\textbf{C}_n}$, the arrow (a natural transformation)
$\vartriangle^1_n(g)$ satisfies the diagonal property 2 of Definition \ref{def:genDiag}, such that for all indexed objects $a_i$ in $\textbf{J}$,

$\vartriangle^1_n(g)(a_j) = g$.\\
Let us show that these diagonal properties are valid also for the comma-propagated functor $\vartriangle_{n+1}:\textbf{C}_{n+1}\rightarrow \textbf{C}_{n+1}^{\textbf{J}}$, and consider an arrow $(g_1;g_2):J(k_1)\rightarrow J(k_2)$ in $\textbf{C}_{n+1} = (\textbf{C}_{n}\downarrow\textbf{C}_{n})$, for which there exists the commutative diagram $g_2\circ k_1 = k_2\circ g_1$ in the lower n-dimensional level $\textbf{C}_{n}$. From the fact that the base functor $\vartriangle:\textbf{C}\rightarrow \textbf{C}^{\textbf{J}}$ is an ascending case, we obtain the object (functor)

$F = \vartriangle_{n+1}(J(k_1)) = L\circ\widehat{\vartriangle}_{n}(J(k_1))$

$= L^0(J\vartriangle^1_{n}\psi)(J(k_1))$

$= L^0J(\vartriangle^1_{n}(k_1))$\\
where $h= \vartriangle^1_{n}(k_1))$ is an arrow (a natural transformation) in $\textbf{C}_n^{\textbf{J}}$. Thus,  from (\ref{eq:L-object}), for each indexed object $a_j$ in $\textbf{J}$,

$\vartriangle_{n+1}(J(k_1))(a_j) = F(a_j) = J(h(a_j))= J(\vartriangle^1_{n}(k_1)(a_j))$

$= J(k_1),~~~~~$ by inductive hypothesis and (\ref{eq:genDiag1}),\\
so that object-component of the coma-propagated functor $\vartriangle_{n+1}$ satisfies the diagonal property. Let us show that the arrow-component of the functor $\vartriangle_{n+1}$ satisfies the diagonal property, that is, for each arrow $(g_1;g_2)$ in $\textbf{C}_{n+1}$, the arrow $\eta = \vartriangle_{n+1}(g_1;g_2)$ in $\textbf{C}_{n+1}^{\textbf{J}}$ is a natural transformation, and hence for each index object $a_j$ in $\textbf{J}$, we have that

$(\vartriangle_{n+1}(g_1;g_2))(a_j)= (L\circ \widehat{\vartriangle}_{n}(g_1;g_2))(a_j)$

$= (L^1(\vartriangle^1_{n}F_{st}^1;\vartriangle^1_{n}S_{nd}^1)(g_1;g_2))(a_j)$

$= (L^1(\vartriangle^1_{n}(g_1);\vartriangle^1_{n}(g_2))(a_j)~~~~~$ for nat.transf. $\vartriangle^1_{n}(g_i), i=1,2$ (arrows in $\textbf{C}_{n}^{\textbf{J}})$

$= (\vartriangle^1_{n}(g_1)(a_j);\vartriangle^1_{n}(g_2)(a_j))~~~~~$  from (\ref{eq:L-object2})

$= (g_1;g_2),~~~~~$ from inductive hypothesis,\\
so that  the arrow-component of the comma-propagated functor $\vartriangle_{n+1}$ satisfies the diagonal property as well.
\\$\square$
\end{example}
 \begin{definition}\label{def:commaPropagationNatTr}
 We define the following four cases of the comma-propagation of a given base natural transformation $\tau:F\rTo^\centerdot G$  (considered as 1-dimensional level denoted by $\tau_1$) which for each (n-1)-dimensional level, $n\geq 2$, will be denoted by $\tau_{n-1}:F_{n-1} \rTo^\centerdot G_{n-1}$ comma-propagated natural transformation, as follows:
 \begin{enumerate}
   \item Basic case: when the base functors are $F,G:\textbf{D}\rightarrow \textbf{C}$. Then n-dimensional  comma-propagated natural transformation is its comma lifted (from (\ref{eq:commaNatTr})),\\
       $\tau_n = \{\tau_{n-1}\}\triangleq~\widehat{\tau}_{n-1}:\widehat{F}_{n-1} \rTo^\centerdot \widehat{G}_{n-1}$, \\that is the function
       $\tau_n  = (\tau_{n-1}F_{st}^0;\tau_{n-1}S_{nd}^0)$, as in the comma-induction.
   \item Descending case: when the base functors are $F,G:\textbf{D}^{\textbf{J}}\rightarrow \textbf{C}$. Then n-dimensional  comma-propagated natural transformation, for the functors in point 2 of Definition \ref{def:commaPropagationFunctor}, is,\\
       $\tau_n = \{\tau_{n-1}\} \triangleq~\widehat{\tau}_{n-1}K:\widehat{F}_{n-1}\circ K \rTo^\centerdot \widehat{G}_{n-1}\circ K$, \\that is, the function
       $\tau_n = ~ (\tau_{n-1}F_{st}^0;\tau_{n-1}S_{nd}^0)K^0$.
   \item Ascending case: when the base functor is $F,G:\textbf{D}\rightarrow \textbf{C}^{\textbf{J}}$.  Then n-dimensional  comma-propagated natural transformation,  for the functors in point 3 of Definition \ref{def:commaPropagationFunctor}, is,\\
       $\tau_n = \{\tau_{n-1}\} \triangleq~L\widehat{\tau}_{n-1}:L \circ\widehat{F}_{n-1}\rTo^\centerdot L\circ\widehat{G}_{n-1}$, \\that is, the function
       $\tau_n =~ L^1(\tau_{n-1}F_{st}^0;\tau_{n-1}S_{nd}^0)$.
   \item Balanced case: when the base functors are $F,G:\textbf{D}^{\textbf{J}}\rightarrow \textbf{C}^{\mathfrak{J}}$  Then n-dimensional  comma-propagated natural transformation,  for the functors in point 4 of Definition \ref{def:commaPropagationFunctor}, is,\\
       $\tau_n = \{\tau_{n-1}\} \triangleq~L\circ \widehat{\tau}_{n-1}K:L \circ\widehat{F}_{n-1}\circ K\rTo^\centerdot L\circ\widehat{G}_{n-1}\circ K$, \\that is, the function
       $\tau_n =~ L^1(\tau_{n-1}F_{st}^0;\tau_{n-1}S_{nd}^0)K^0$.
 \end{enumerate}
 \end{definition}
So, as in the case of the comma-propagated functors, also for the comma-propagated natural transformations different from basic case, we have the modulation by the functors $K$ and $L$.
 It is easy to verify that in all four cases above,
  we obtain that for each $n\geq 1$, the comma-propagated natural transformation is just between comma-propagated functors derived from the base functors $F$ and $G$:

 $\tau_n:F_n \rTo^\centerdot G_n$\\
 The following distributive laws \cite{Majk23s} between the operation of comma-propagation '$\{\_\}$'  and the following categorial operations, for $n\geq 1$:
 \begin{enumerate}
   \item The operation of composition of functors:\\
   $F_n\circ G_n  = (F\circ G)_n$
   \item The horizontal composition of natural transformations:\\
   $\tau_n\circ \eta_n =(\tau\circ \eta)_n$
   \item The vertical composition of natural transformations:\\
   $\tau_n\bullet \eta_n = (\tau\bullet \eta)_n$
 \end{enumerate}
\section{Categorial Topology and Symmetries under Comma-propagation Transformations}
 In mathematics, topology  is concerned with the properties of a geometric object that are preserved under continuous deformations, such as stretching, twisting, crumpling, and bending; that is, without closing holes, opening holes, tearing, gluing, or passing through itself.  A property that is invariant under such deformations is a topological property.
 The motivating insight behind topology is that some geometric problems depend not on the exact shape of the objects involved, but rather on the way they are put together.
 In one of the first papers in topology, Leonhard Euler demonstrated that it was impossible to find a route through the town of K\"{o}nigsberg (now Kaliningrad) that would cross each of its seven bridges exactly once. This result did not depend on the lengths of the bridges or on their distance from one another, but only on connectivity properties: which bridges connect to which islands or riverbanks. This Seven Bridges of K\"{o}nigsberg problem led to the branch of mathematics known as graph theory.
 But, the graph theory is fundamental part of the category theory, that is, each graph can be extended into a category as, for example the small index categories $\textbf{J}$ derived from special graphs and used for co(limits) and categories of functors $\textbf{C}^{\textbf{J}}$.

 Consequently, in Categorial Topology, given a category (as a "geometric object") $\textbf{C}$ we can consider its properties preserved under continuous action (a "deformation") of a comma-propagation operation $\{\_\}$ provided in previous Section. Indeed, we can consider any category $\textbf{C}$ as an abstract geometric object \cite{Majk24G}, that is a discrete space where the points of such abstract space are the objects of this category and arrows between objects as the paths: given any two points (two objects of the category) we can have a number of oriented paths from first to the second point.some of them equal (commutative diagrams in this category between these two objects).

 In geometry, groups arise naturally in the study of symmetries and geometric transformations: The symmetries of an object form a group, called the symmetry group of the object, and the transformations of a given type form a general group. Lie groups appear in symmetry groups in geometry, and also in the Standard Model of particle physics. The Poincare group is a Lie group consisting of the symmetries of spacetime in special relativity. Point groups describe symmetry in molecular chemistry.

 By consecutive application of the comma-propagation operation $\{\_\}$ in Definitions \ref{def:commaPropagationFunctor} and \ref{def:commaPropagationNatTr}, we can obtain a number of infinite sets $X$:
 \begin{definition}\label{def:CatSymX}
  Let $X$ be an infinite set of all n-dimensional categories $\textbf{C}_n$, or of categories of functors $\textbf{C}_n^{\textbf{J}}$, or of functors $F_n$ or of natural transformations $\tau_n$, $n\geq 1$, so that comma-propagation can be defined as the function $\{\_\}:X \rightarrow X$.\\
   Then we can define the following subsets of $X = \{a_1,a_2,a_3,...\}$, for each $k\geq 1$,
  \begin{equation}\label{eq:CatSymSubset}
 X_k \triangleq \{a_m \mid a_m \in X, m\geq k+1\}
 \end{equation}
  obtaining the following chain of infinite subsets of $X$,

  $X \supset X_1 \supset X_2 \supset X_3 \supset ....$
 \end{definition}
 Let us show how also in Category Theory the concept of the symmetries is related to the group of transformations derived from the comma-propagation and what are more significant invariances under such transformations.
\begin{definition}\label{def:CatSymGroup}
 We define the following infinite abelian category-symmetry group $CS(\mathbb{Z}) =(\cdot, g_0,g_1, g_{-1}, ...)$ of comma-transformations, with composition operation $\cdot$, identity element $g_0$, comma-up transformation element $g_1$ (whose action on a given set $X$  in Definition \ref{def:CatSymX} is the comma-propagation operation (function) $\widetilde{g}_1 = \{\_\};X \rightarrow X$), its inverse comma-down transformation element $g_{-1}$, and for each $k\geq 2$,
 \begin{equation}\label{eq:CatSymGroup}
 g_k \triangleq \overbrace{g_1\cdot ...\cdot g_1}^k   ~~~and ~~its~~inverse~~~ g_{-k} \triangleq \overbrace{g_{-1}\cdot ...\cdot g_{-1}}^k
 \end{equation}
 such that for any two elements $g_n$ and $g_m$, $n,m \in \mathbb{Z}$, we have that $g_k\cdot g_m = g_{n+m}$.
\end{definition}
So, if we denote by $(\mathbb{Z},+) =(+,0,1,-1,...)$ the group of integers with addition  as group operation and $0$ as identity element, then we have simple isomorphism of groups:
\begin{equation}\label{eq:CatSymGroupIso}
 \sigma;(\mathbb{Z},+)\simeq CS(\mathbb{Z})
 \end{equation}
 Let us consider now the actions of the category-symmetry group on the sets of  n-dimensional categorial concepts specified in Definition \ref{def:CatSymX}:
 \begin{definition}\label{def:CatSymAction}
 We define the action of the $CS(\mathbb{Z})$ group on any infinite set $X$ of n-dimensional categorial concepts specified in Definition \ref{def:CatSymX}, as possibly partial function $CS(\mathbb{Z}) \times X \rightarrow X$, written, for any $n,m \in \mathbb{Z}$ and $k\geq 1$, as

 $(g_n,a_k) \mapsto g_na_k~~~$ if $k+n \geq 0$,\\
 so that $g_0a_k =a_k$, but generally $(g_n\cdot g_m)a_k \neq g_n(g_m a_k)$ because if $m+k < 0$ then $g_m a_k$ is not well defined.
 Generally, for each group element $g_n \in CS(\mathbb{Z})$, if $n\geq 0$ we have the (total) function

 $\widetilde{g}_n:X\rightarrow X~~~$ such that for each $a_m \in X$ with $k\geq 1$, $\widetilde{g}_n(a_m) = a_{n+m}$\\
 while for each negative $n< 0$, we obtain a partial function on $X$,

 $\widetilde{g}_n:X_{-n}\rightarrow X~~~$ such that for each $a_m \in X_{-n}$, i.e, $n+m \geq 1$,  $\widetilde{g}_n(a_m) = a_{n+m}$\\
 where the subset $X_{-n} \subset X$ is given by (\ref{eq:CatSymSubset}).
 \end{definition}
 This partiality of the action of the category-symmetry group to n-dimensional categorial concepts, with $n\geq 1$ only, can be avoided reducing our attention only to comma-up transformations from initial (base) concepts  $\textbf{C}_1$, or of categories of functors $\textbf{C}_1^{\textbf{J}}$, or of functors $F_1$ or of natural transformations $\tau_1$, into higher n-dimensional concepts in order to verify what categorial properties of the basic category $\textbf{C} = \textbf{C}_1$ preserves (remains invariant) under the comma-propagation transformations in all higher n-dimensional levels $\textbf{C}_n$ for $n = 2,3,4,...$.

 So, our categorial symmetry can be described by the Peano-like system of axiomatic definition of natural numbers $n = 1,2,3,...$, where instead of a number $n$ we have the n-dimensional level $\textbf{C}_n$ and instead of Peano successor operation "+1" the comma-propagation operation $\{\_\}$, with the following Peano-like axioms for n-dimensional levels:
 \begin{enumerate}
   \item $\textbf{C}_1 = \textbf{C}$ is a n-dimensional level.
   \item Every n-dimensional level $\textbf{C}_n$ has a successor which is also a n-dimensional level  $\textbf{C}_{n+1} = \{\textbf{C}_n\}$.
   \item $\textbf{C}_1$ is not successor of any n-dimensional level.
   \item If the successor of n-dimensional level $x$ equals to the successor of n-dimensional level $y$, then $x$ equals $y$.
    \end{enumerate}
 Thus, the principle of induction can be used to determine the fundamental invariance  of categorial statements (which generally is represented by kinds of universal categorial commutative diagrams (that is, by mathematical equations of the composition of arrows in a given category): If a category statement is true for $\textbf{C}_1$ (in base category $\textbf{C}$), and the truth of this statement implies its truth for the successor $\textbf{C}_2 = \{\textbf{C}_1\}$, then the statement is true for every n-dimensional level $\textbf{C}_n$, $n\geq 2$.

 In this, Peano-like approach to fundamental categorial symmetry based on comma-propagation transformations in the space of n-dimensional levels, instead of the symmetry group $CS(\mathbb{Z})$, we can use only the categorial-symmetry commutative monoid $CS(\mathbb{N}) =(\cdot, g_0,g_1,g_2,....)$ (which is a closed subalgebra of $CS(\mathbb{Z})$), that is, its semigroup with the same identity element $g_0$), where  $\mathbb{N} = \{0,1,2,...\} \subset \mathbb{Z}$ is the set of natural numbers, with the following isomorphism of  monoids:
 \begin{equation}\label{eq:CatSymGroupIso}
 \sigma_0;(\mathbb{N},+)\simeq CS(\mathbb{N})
 \end{equation}
 Let us consider now the actions of the category-symmetry group on the sets of  n-dimensional categorial concepts specified in Definition \ref{def:CatSymX}:
 \begin{definition}\label{def:CatSymActionM}
 We define the action of the categorial-symmetry  $CS(\mathbb{N})$ monoid on any infinite set $X$ of n-dimensional categorial concepts specified in Definition \ref{def:CatSymX}, as a function $CS(\mathbb{N}) \times X \rightarrow X$, written, for any $g_n,g_m \in CS(\mathbb{N})$ and $a_k\in X$, as

 $(g_n,a_k) \mapsto g_na_k~~~$ \\
 so that $g_0a_k =a_k$,  and $(g_n\cdot g_m)a_k = g_n(g_m a_k)$.
 \end{definition}
  In this way, as shown in the book \cite{Majk23s}, by the action of the symmetry group $CS(\mathbb{Z})$ or symmetry monoid $CS(\mathbb{N})$ on the sets $X$ of the n-dimensional categorial concepts specified in Definition \ref{def:CatSymX}, we can show the validity of the Peano-like principle of induction for categorial statements as co(limits) and categorial adjunctions, and their invariance under comma-propagation transformations.

Adjunctions, $(F,G,\varepsilon,\eta)$ of the functor $F:\textbf{C}\rightarrow \textbf{D}$ with its right-adjoint functor $G:\textbf{D}\rightarrow \textbf{C}$ which generate two natural transformations (the unit $\eta:I_\textbf{C}\rTo^\centerdot GF$ and counit $\varepsilon:FG \rTo^\centerdot I_\textbf{D}$, which in the dual version of the adjunction have opposite arrow direction) introduced by Kan in 1958, provide a descriptive framework of great generality, capturing the essence of many canonical constructions.
They turn up throughout mathematics often as "closures" and "completions", and as "free" and "generated" structures. For example, the transitive closure of a graph, the completion of a metric space, factor commutator groups, and free algebras are all examples of adjunctions. In categorical logic, quantifiers are interpreted as adjunctions with respect to substitution of variables.

The canonical nature of these constructions is captured by universality. Indeed, adjunctions subsume the universal structures, both limits and colimits, \index{(co)limits} when  $\textbf{D} = \textbf{C}^{\textbf{J}}$ is a category of functors for a small indexed category $\textbf{J}$ and $F= \bigtriangleup$ is the diagonal functor. Fundamentally, the adjunctions represents two "forces" generated by the composition of these two functors creating the unit natural transformation such that for each object $c$ in the category $\textbf{C}$ is generated the displacement (arrow) toward the object $GF(c)$, and by the counit natural transformation such that for each object $d$ in the category $\textbf{D}$, is generated the displacement (arrow) from the object $FG(d)$ into this object $d$. From the fact that these forces are applied to \emph{all} objects of these two categories, we can consider an adjunction as a kind of the scalar \emph{fields} in these two categories: \emph{adjunctions-as-fields}. \index{adjunctions-as-fields} More about this interpretation of the adjunction has been discussed in \cite{Majk23s}.

(Co)limits are a generalization of the universally defined structures. The generalization is based upon diagrams in a category.
We define diagrams and then the general structure of a (co)limit. Limits and colimits can both be defined through universality. More
examples will be found in next. A curious result is the
interdefinability of universal concepts in category theory. Universal concepts
can be obtained from one another by interpreting them in suitable
categories. Limits and colimits give a unified treatment of constructs such as
products and sums of pairs of objects and limits of chains of objects.

 Both of these important categorial concepts are built by using the more simple concepts of universal arrows. By considering that, from Theorem \ref{th:A1}, the comma lifting preserves the universal arrows, the next theorem (proof can be find in \cite{Majk23s}) demonstrates the universal arrows are preserved by the comma-propagation: \index{comma-propagation transformation}
 \begin{theo} \label{th:A7}
 Given a functor $G$,
 for each universal arrow "\emph{from $c$ to functor $G$}" which is the pair $(d,g)$, with object $d$ in domain category of $G$ and arrow $g:c\rightarrow G(d)$ in the codomain category of $G$, presented by the figure\footnote{Where for each pair $(d',f)$ there exists \emph{the unique} arrow $\underline{f}:d\rightarrow d'$ such that $f = G(\underline{f})\circ g$.},
  \begin{equation} \label{fig:universal-arrow2n}
\begin{diagram}
 c         &\rTo^{g}          & G(d)    && d\\
           & \rdTo_{f}&     \dTo_{G(\underline{f})}&&\dDashto_{\underline{f}}     \\
           &                 &   G(d')   &&d'
\end{diagram}
\end{equation}
with $G_1 = G$, there is the comma-propagated universal arrows as follows:
 \begin{enumerate}
   \item Basic case: when $G:\textbf{D}\rightarrow \textbf{C}$ then we have the universal arrow $(J(id_d), (g;g))$ from $J(id_c)$ to comma-propagated  functor $G_2:(\textbf{D}\downarrow\textbf{D})\rightarrow (\textbf{C}\downarrow\textbf{C})$.
   \item Descending case: when $G:\textbf{D}^{\textbf{J}}\rightarrow \textbf{C}$ then we have the universal arrow $(L^0J(id_d), (g;g))$ from $L^0J(id_c)$ to comma-propagated  functor $G_2:(\textbf{D}\downarrow\textbf{D})^{\textbf{J}}\rightarrow (\textbf{C}\downarrow\textbf{C})$.
   \item Ascending case: when $G:\textbf{D}\rightarrow \textbf{C}^{\textbf{J}}$ then we have the universal arrow $(J(id_d), L^1(g;g))$ from $J(id_c)$ to comma-propagated  functor $G_2:(\textbf{D}\downarrow\textbf{D})\rightarrow (\textbf{C}\downarrow\textbf{C})^{\textbf{J}}$.
   \item Balanced case: when $G:\textbf{D}^{\textbf{J}}\rightarrow \textbf{C}^{\mathfrak{J}}$ then we have the universal arrow $(L^0J(id_d), L^1(g;g))$ from $L^0J(id_c)$ to comma-propagated  functor $G_2:(\textbf{D}\downarrow\textbf{D})^{\textbf{J}}\rightarrow (\textbf{C}\downarrow\textbf{C})^{\mathfrak{J}}$.
 \end{enumerate}
 where $L$ is modulator functor isomorphism specified in point 1 of Definition \ref{def.K-L-iso}.

  These results hold also for the dual concept of co-universal arrow "\emph{from functor $G$ to $c$}",
  with opposite direction of all arrows in figure (\ref{fig:universal-arrow2n}) above.
 \end{theo}
For example, let us consider the following co-universal arrows, for general diagonal functor in Definition \ref{def:genDiag} with $\textbf{C}= \textbf{D}\times \textbf{D}$  and $\textbf{J} = \textbf{1}$ is the category of only one object $a_1$ and its identity arrow $id_{a_1}$ with figure (\ref{fig:DiagGen}).
 \begin{example} \label{ex:CommaPropProduct}
 In this case $\vartriangle:\textbf{D}\rightarrow (\textbf{D}\times \textbf{D})^{\textbf{1}}$ (which is an ascending functor, as discussed in details in Example (\ref{exam:DiagFunct})).

  Let us consider the opposite descending functor $G:(\textbf{D}\times \textbf{D})^{\textbf{1}}\rightarrow \textbf{D}$ with \emph{dual concept} of co-universal arrow "\emph{from functor $G$ to $c$}" (with opposite arrows w.r.t figure in
(\ref{fig:universal-arrow2n})), represented by the pair $(d,g)$ such that

 $d(a_1) \triangleq (c,c)$,\\
  represented by  the figure below, with its  re-elaborated part here and with

   $d'(a_1) \triangleq (a,b)$,
\begin{equation} \label{fig:Cartesian productUnivArrow2}
\begin{diagram}
 c &\lTo^{
 g}              & G(d) &&& &      d:\textbf{1}\rightarrow \textbf{D}\times \textbf{D}&& &    (c,c) = d(a_1)   \\
     & \luTo_{f}&  \uTo_{G(\underline{f})} &&& &\uDashto^\centerdot_{\underline{f}}&  &&    \uTo^{(k_1,k_2)}_{=\underline{f}(a_1)}  \\
     & &    G(d')   &&&& d':\textbf{1}\rightarrow \textbf{D}\times \textbf{D}&&&        (a,b) =d'(a_1) \\
      & & \textbf{D}&&\lTo^G&&   (\textbf{D}\times \textbf{D})^{\textbf{1}}   & &&  \textbf{D}\times \textbf{D}
\end{diagram}
\end{equation}
In the case when $G$ is the coproduct functor, so that

$G(d) = c+c$,

$G(d') = a+b$,

 $G(\underline{f}) = k_1+k_2$ \\
  so that

 $g = [id_c,id_c]$\\
  and hence

  $f = g\circ G(\underline{f}) = [k_1,k_2]$, \\
   we obtain an alternative definition of the commutative diagram  
   representing the coproducts by the adjunction of the alternative (more general) diagonal functor $\vartriangle:\textbf{D}\rightarrow \textbf{D}^{\textbf{2}_D}$ where $\textbf{2}_D$ is the discrete index category composed by two index-objects $a_1$ and $a_2$, and the right adjoint to it modified functor $~G:\textbf{D}^{\textbf{2}_D} \rightarrow \textbf{D}$.

In this case for the coproducts we can use also the standard definition by functor $+:\textbf{D}\times \textbf{D}\rightarrow \textbf{D}$, from the simple product category $\textbf{D}\times \textbf{D}$ as presented in diagram above.

With this we explained how the standard diagonal functor $\vartriangle:\textbf{D}\rightarrow \textbf{D}\times \textbf{D}$ can be replaced in the large family of general diagonal functors $\vartriangle:\textbf{D}\rightarrow (\textbf{D}\times \textbf{D})^{\textbf{J}}$, as we explained it by figure (\ref{fig:arrow-diagonal}) for the representation of the  arrow-diagonal functor $\blacktriangle:\textbf{D}\rightarrow (\textbf{D}\downarrow \textbf{D})$, and explains the mathematical importance for introduction of the more powerful general diagonal functors.
\\$\square$
 \end{example}
Advantage of using the general diagonal functors $\vartriangle:\textbf{D}\rightarrow \textbf{C}^{\textbf{J}}$ is  that they are ascending cases in the comma-propagation of functors as well, and hence they are a consistent part of the comma-propagation theory, important for the extension of the comma-propagation to all (co)limits (provided in next section).

 The properties demonstrated previously  for the comma-propagation of functors, natural transformations and (co)universal arrows within the infinite hierarchy of n-dimensional levels (arrow categories), are necessary in order to expand the comma-propagation theory to the most useful category structures: adjunctions and (co)limits.  Category theory is the theory of typed composition of arrows and as such is a very "weak" theory. The theory finds its strength in powerful and intricate descriptive mechanisms.

Thus, from the fundamental role of adjunctions and (co)limits,  we will show how these important categorial structures, presented in the base categories, propagates by comma-propagation to all infinite hierarchy of n-dimensional levels (the proof is provided in \cite{Majk23s}).
%
\begin{propo} \label{prop:Comma-PropAdjunct} \textsc{Comma-propagation of adjunctions}:\\
Let $\textbf{C}^{\textbf{J}}$  be a category of functors from a small index category $\textbf{J}$ into a base category $\textbf{C}$. Then, for every adjunction $(F,G,\varepsilon,\eta)$ where $F:\textbf{D}\rightarrow \textbf{C}^{\textbf{J}}$ is left adjoint of a functor $G$  with $\eta$ and $\varepsilon$ the unit and counit of this adjunction,  for each n-dimensional level $\textbf{C}_n$, $n\geq 1$, there exists an adjunction between comma-propagated functors, $(F_n,G_n,\varepsilon_n,\eta_n):\textbf{D}_n\rightarrow \textbf{C}_n^{\textbf{J}}$.
\end{propo}
 Based on this result of comma-propagation of the adjunctions, and on the result shown in the Example \ref{exam:DiagFunct}  that the comma-propagation of the general diagonal functor $\vartriangle:\textbf{C}\rightarrow \textbf{C}^{\textbf{J}}$ preserves its diagonal functor properties, and from the fact that it will be used for consideration of the (co)limits in the n-dimensional hierarchy, we can show this comma-propagation of the (co)limit structures:
\begin{theo}\label{th:A.8}
Let $\textbf{C}^{\textbf{J}}$ be a category from a small index category $\textbf{J}$ into a given base category  $\textbf{C}$, such that for any object (functor) $d':\textbf{J}\rightarrow \textbf{C}$ in $\textbf{C}^{\textbf{J}}$, there exists a limit $G(d')$ in $\textbf{C}$ generated by the functor $G:\textbf{C}^{\textbf{J}} \rightarrow \textbf{C}$ right adjoint to the general diagonal functor $\vartriangle:\textbf{C}\rightarrow \textbf{C}^{\textbf{J}}$.

Consequently, for every functor $d':\textbf{J}\rightarrow \textbf{C}_n$ which is an object in the category of functors $\textbf{C}_n^{\textbf{J}}$, for $n\geq 2$, there exists a limit $G_n(d')$ in  $\textbf{C}_n$, where $G_n:\textbf{C}_n^{\textbf{J}}\rightarrow  \textbf{C}_n $ is the comma-propagated functor  of  $G$ and right-adjoint to the comma-propagated diagonal functor $\vartriangle_n:\textbf{C}_n\rightarrow \textbf{C}_n^{\textbf{J}}$.
\end{theo}
\textbf{Proof}: By using the distributive laws at the end of previous section and  Proposition
\ref{prop:Comma-PropAdjunct} with adjunction $(\vartriangle,G,\varepsilon,\eta)$ where $\vartriangle:\textbf{D}\rightarrow \textbf{C}^{\textbf{J}}$ is general diagonal functor left adjoint of a functor $G$ its right adjoint functor $G$ with $\eta$ and $\varepsilon$ the unit and counit of this adjunction, so that, \emph{for any object} $c$ in $\textbf{C}$, we have the universal arrow, i.e., pair $(d,g)$ where $d =(d^0,d^1)=\vartriangle^0(c)$ is the constant functor from a small index category $\textbf{J} = \textbf{Sch}(graph)$  (derived as a sketch category from a finite  graph, which in this case represent the structure of the diagram of the limit cones) into $\textbf{C}$ and $g=\eta(c)$,  such that for each pair $(d',f)$ there is a unique arrow (natural transformation) $\underline{f}:d\rTo^{\centerdot} d'$ such that the following diagram commutes
\begin{equation} \label{fig:universal-arrow40CP}
\begin{diagram}
 c~         &\rTo^{\eta(c)}          &~ G(d)    &&& &d = \vartriangle(c)&\rTo^{\vartriangle(f)}_\centerdot&\vartriangle G(d') \\
           & \rdTo_{f}&     \dTo_{G(\underline{f})}&&&&\dTo_\centerdot^{\underline{f}}&\ldTo^\centerdot_{\varepsilon(d')} &     \\
           &                 &   G(d')   &&&&d':\textbf{J}\rightarrow \textbf{C}\\
            & & \textbf{C} &&\lTo^{G}&& \textbf{C}^{\textbf{J}}
\end{diagram}
\end{equation}
where the object (functor) $d':\textbf{J}\rightarrow \textbf{C}$ in $\textbf{C}^{\textbf{J}}$ represents the diagram for which the object $G(d')$ in $\textbf{C}$ is the limit (that is, \emph{$G$ creates the limit of the diagram $d'$}).

From the fact that for each index object $a_j$ in $\textbf{J} = \textbf{Sch}(graph)$, we have for the constant functor $d = \vartriangle(c)$, from (\ref{eq:genDiag0}), $d^0(a_j)= c$ and each arrow $l_{jm}$ in $\textbf{J}$, $d^1(l_{jm}) = id_c$  from Definition \ref{def:genDiag}) of general diagonal functor,
it holds that the diagram represented by this functor $d = \vartriangle(c)$ is composed by all object equal to $c$ and all arrows equal to its identity arrow $id_c$, so that the diagram represented by this functor $d = \vartriangle(c)$  can be reduced to the single object $c$ which is a vertex of the cone composed by arrow components $\underline{f}(a_j):c\rightarrow d'(a_j)$ for each index object $a_j$ in  $\textbf{J}$, where $d'(a_j)$  is an object of the diagram of the graph represented by the functor $d':\textbf{J}\rightarrow \textbf{C}$.

That is,  for each object $c$ in $\textbf{C}$, the  arrow (a natural transformation) $\underline{f}:d \rTo^\centerdot  d'$ represents a cone with vertex equal to object $c$  toward the diagram represented by the functor $d'$ (the arrow components $\underline{f}(a_j):c \rightarrow d'(a_j)$ for each index-object $a_j$ in $\textbf{J}$ generate this cone). The \emph{limit cone} is represented by the natural transformation $\varepsilon (d')$ derived from the counit $\varepsilon$ for thie given diagram (a functor) $d'$,  with the vertex equal to the limit $G(d')$ of this diagram $d'$ (from the fact that the functor $\vartriangle (G(d'))$ reduces to the single object $G(d')$ in $\textbf{C}$ as explained previously for the diagonal functor $\vartriangle:\textbf{C} \rightarrow \textbf{C}^{\textbf{J}}$.

From the fact that the comma-propagation preserves the adjunctions, the diagram above holds for each comma-propagated adjunction $(\vartriangle_n,G_n,\varepsilon_n,\eta_n):\textbf{C}_n\rightarrow \textbf{C}_n^\textbf{J}$, for $n\geq 1$, it holds that the limits in $\textbf{C}$ are comma-propagated to all higher n-dimensional levels $\textbf{C}_n$.

The proof for the dual case of the \emph{colimits} is analog to this done above, by considering that in this way all arrows in the commutative diagrams in (\ref{fig:universal-arrow40CP}) \emph{are inverted}.
\\$\square$\\
Finally, we obtain the following closure property for the n-dimensional levels:
\begin{coro} \label{coro"small-complete}
Let $\textbf{C}$ be a small-complete category. Then, every n-dimensional level $\textbf{C}_n$, for $n\geq 2$, is a small-complete category as well.
\end{coro}
\textbf{Proof}: Directly from Theorem \ref{th:A.8}, by considering  every small index category (derived from a finite diagram) $\textbf{J}$.
\\$\square$\\
It is well known that a category $\textbf{C}$ is a small-complete if every its finite diagram has  a limit: the enough condition that it is a small complete is that $\textbf{C}$ has arbitrary finite products and every pair of its arrows has an equalizer.

Thus, these two particular limits, the products and the equalizers have particular importance. One of the well known categories with such properties is the \textbf{Set} category.
\section{Conclusions \label{sec:ExamplCP}}
We provided an abstract category theory for the arrow categories and \emph{global symmetries} of comma-propagation transformations with an introduction of the infinite hierarchy  of the arrow categories considered as "\emph{Galilean reference frames}" of a given base category $\textbf{C}$ (the n-dimensional levels $\textbf{C}_n$, for $n \geq 1$), the two fundamental comma projection-functors with natural transformation $\psi$ between them and inverse to it encapsulation operator $J$. We provides a new categorial-based (non set-based) definition of natural numbers.  Then we introduce the basic transformations of functors and natural   along this abstract dimension, from n-dimensional to (n+1)-dimensional level, denominated as "comma lifting".

Here we extended this approach to the more general phenomena of \emph{comma-propagation transformations} (a "Galilean boost" of "uniform motion" between two n-dimensional levels seen as different reference frames) and its invariant properties. Comma-propagation is a general transformation between n-dimensional levels, applied to the functors and natural transformations. If we consider the universal properties of a base category  (universal arrows, adjunctions) as a kind of different "Lagrangian" in the Category Theory, then we can consider their invariance under this comma-propagation transformations: a \emph{global symmetry} in this case means that these categorial universal structures are preserved under these general transformations \emph{in all} n-dimensional levels (the law of conservation). The application of this kind of symmetry are provided as well in dedicated sections in \cite{Majk23s}.

%

%

%

\end{document}

%% file: macros-general.tex
\begingroup
\catcode`\~=11
\gdef\urltilde{\lower 0.6ex\hbox{~}}
\endgroup

 \renewcommand{\P}{\mathcal{P}}

